\documentclass[11pt]{article}

\usepackage[cp1251]{inputenc}
\usepackage[russian]{babel}
\usepackage{amsthm,amsfonts,amsbsy, amssymb,amsmath,graphicx,euscript}
\usepackage{makeidx}
\usepackage{multicol}

\usepackage{mathrsfs}

\usepackage{amsfonts}
\usepackage{amssymb}
\usepackage{amsmath}
\usepackage{epsfig}
\usepackage{psfrag}
\usepackage{floatflt}
\usepackage{mathrsfs}
\author{I. D. Remizov, A. V. Savvateev}

\newcommand{\Argmax}{\mathrm{Argmax}}

\begin{document}

\begin{center}\Large\textbf{D(Maximum) = P(Argmaximum)}

\end{center}

\begin{center}

\textbf{Ivan D. Remizov}\\ 

\footnotesize

Department of Function Theory and Functional Analysis, Faculty of Mechanics and Mathematics,\\ Moscow State University, Russia\\
Email: \texttt{ivremizov@yandex.ru} $\ $  Homepage: \texttt{http://ivanremizov.ru}

\normalsize

\end{center}

\begin{center}

\textbf{Alexei V. Savvateev}\\

\footnotesize

New Economic School, Moscow, Russia\\
Email: \texttt{hibiny@mail.ru} $\ $  Homepage: \texttt{www.nes.ru/\~{}savvateev}

\end{center}
\normalsize

\footnotesize

\textbf{Abstract.} In this note, we propose a formula for the subdifferential of the maximum functional $m(f) = \max_K f$ on the space of real-valued continuous functions $f$ defined on an arbitrary metric compact $K$. We show that, given $f$, the subdifferential of $m(f)$ always coincides with the set of all probability measures on the arg-maximum (the set of all points in $K$ at which $f$ reaches the maximal value). In fact, this relation lies in the core of several important identities in microeconomics, such as Roy's identity, Sheppard's lemma, as well as duality theory in production and linear programming.\\

\normalsize

Let $K$ be a compact metric space, and let $L=C(K)$ be the vector space of all continuous real-valued functions on $K$ endowed with the uniform norm $\|x\|=\sup_{t\in K}|x(t)|$. Let $L'$ be the vector space of all continuous real-valued linear functionals on $L$ endowed with the operator norm $\|h\|=\sup_{x\in L, x\neq 0}\frac{|h(x)|}{\|x\|}$. 

Let $m$ be the (convex) functional on $L$ that assigns to each function $f\colon K\to\mathbb{R}$ its' maximal value $m(f)=\max_{t\in K} f(t)$. 

One can define the sub-gradient of $m$ at a point $f$ as $$\partial m_f=\left\{h\in L':m(g) - m(f)\geq  h(g-f)\ \forall g\in L\right\}\subset L'.\eqno(1)$$

Let us denote $\Argmax_K f=f^{-1}(m(f))=K_f$. Obviously $K_f$ is a closed (and thus compact) subset of $K$.

It is known that each real-valued sigma-additive measure $\mu$ on the sigma-algebra $\mathcal{A}$ of all Borelian subsets of $K$ defines a real-valued continuous linear functional on $L$ in such a way:$h\colon x\longmapsto \int_K x(t)\mu(dt)$. By the Reisz-Markov theorem (Dunford, Schwartz~1958) there exists an isomorphism between linear normalized spaces $L'$ and $\mathfrak{M}(K)$, where $\mathfrak{M}(K)$ is the space of all real-valued sigma-additive Borelian measures on $K$ endowed with the measure variation norm. 

Let us denote the set of all probability measures on $K$ concentrated on $K_f$ by  $$\mathcal{P}_f=\{\mu\in\mathfrak{M}(K): \mu K_f=1, \mu A\geq 0, \mu A=\mu(A\cap K_f)\ \forall A\in\mathcal{A}\}\subset L'.$$

\textbf{Proposition.} Under these conditions, for each $f\in C(K)$, the sub-gradient 
of the maximum function $\max_K f$ coincides with the set of all probability measures 
on $\Argmax_K f$, i.e.

$$\partial m_f=\mathcal{P}_f\ \forall f\in L.$$
 
\textbf{Proof.} Let us consider $f\in L$ and prove the inclusion $\partial m_f\supset\mathcal{P}_f$. Indeed, let $\mu$ be a probability measure on $K$ concentrated on $K_f$. Let us denote $h=\left[x\longmapsto \int_{K_f} x(t)\mu(dt)\right]\in L'.$ One can see that, since $\mu K_f=1$ and function $f$ is constant and equal $m(f)$ on $K_f$, $m(f)=\int_{K_f}f d\mu$ . 

Consider arbitrary function $g\in L$. Estimate $$m(g)-m(f)\geq \int_{K_f}g\ d\mu - m(f)=\int_{K_f}g\ d\mu - \int_{K_f}f\ d\mu=\int_{K_f}(g-f)\ d\mu=h(g-f)$$ and (1) imply $h\in\partial m_f$.\\

Now let us consider $f\in L$ and prove that $\partial m_f\subset\mathcal{P}_f$. Indeed, by the Reisz-Markov theorem, for each functional $h\in \partial m_f\subset L'$, there exists a real-valued sigma-additive measure $\mu$ on $\mathcal{A}$ such that $h(x)= \int_K x(t)\mu(dt)$ $\forall x\in L$. Let us show that $\mu$ is a probability measure and that it is concentrated on $K_f$.

One can set $g= \mathbf{0}$ in (1), and by the linearity of $h$ conclude that $m(f)\leq h(f)$. For every $x\in L$ one can set $g=f+x$ in (1) and see that $m(f+x)-m(f)\geq h(x)$. Setting $x=f$ leads us to  $m(f)\geq h(f)$. Therefore, $$h(f)=m(f).\eqno(2)$$

We are eager to prove the non-negativity of the mesure $\mu$. It follows from the Reisz-Markov theorem that non-negative functional (i.e. functional which takes non-negative value on any non-negative function) corresponds to a non-negative measure, and vice versa. So let us prove that for each non-negative function $x\in L$ inequality $h(x)\geq 0$ holds. First of all, one can see that non-negativity of $x$ implies $m(f)\geq m(f-x),$ i.e. $m(f) - m(f-x)\geq 0.$ Now let us set $g=f-x$ in (1) and see, that $m(f-x)-m(f)\geq h(-x);$ functional $h$ is linear, thus from the above estimates inequality $h(x)\geq m(f)-m(f-x)\geq 0$ can be easily derived. This proves that measure $\mu$ is non-negative since functional $h$ is non-negative, i.e. proves that $$\mu A\geq 0\ \forall A\in\mathcal{A}.\eqno(3)$$

Let us show that $\mu$ is a probability measure. One can see that the estimate $$h(f)=\int_K f(t)\mu(dt)\stackrel{(3)}{\leq} \int_K \max_{t\in K}f(t)\mu(dt)=\int_K m(f) d\mu=m(f)\mu K\stackrel{(2)}{=}h(f)\mu K$$ implies  $1\leq \mu K$. Setting $g=f+\mathbf{1}$ in (1) leads us to the inequality $1\geq h(\mathbf{1})=\mu K$. Therefore, $$\mu K=1.\eqno(4)$$

It remains to prove that measure $\mu$ is concentrated on $K_f,$ i.e. that $\mu K_f=1.$ It follows from (4) and from the additivity of $\mu$ that $\mu K_f=1$ is equivalent to $\mu(K\setminus K_f)=0$ and to $\mu A=\mu (A\cap K_f)\ \forall A\in\mathcal{A}.$ Relations $$m(f)\stackrel{(2)}{=}h(f)=\int_K f(t)\mu(dt)\stackrel{(3)}{\leq} \int_K m(f)\mu(dt)=m(f)\mu K\stackrel{(4)}{=}m(f)$$ imply that $\int_K (m(f)-f)\ d\mu=0$. Consider the function $g(x)=m(f)-f(x);$ it is equal to zero on $K_f$ and it is positive on $K\setminus K_f$. However $\int_{K\setminus K_f}g\ d\mu=0$, which is possible only if $\mu (K\setminus K_f)=0$.\\ 

The proof is complete.

\section*{References} 

Dunford N., Schwartz J. T. \textit{Linear operators. Part I: general theory.} --- Interscience publishers, New York, London, 1958.
 
\end{document}